\documentclass{amsart}
\usepackage{mathtools}
\usepackage{url}
\usepackage{enumitem}

\usepackage{tikz}

\usepackage{xparse}

\let\origmaketitle\maketitle
\def\maketitle{
  \begingroup
  \def\uppercasenonmath##1{} 
  \let\MakeUppercase\relax 
  \origmaketitle
  \endgroup
}
    \providecommand\given{}%
\newcommand\SetSymbol[1][]{%
    \nonscript\:#1\vert%
    \allowbreak%
    \nonscript\:%
    \mathopen{}}%
\DeclarePairedDelimiterX\Set[1]\{\}{%
    \renewcommand\given{\SetSymbol[\delimsize]}%
    #1%
}%
\newcommand*{\RealSpace}{%
    \mathbb{R}^3%
}%
\newcommand*{\Spaceradius}{r}%
\NewDocumentCommand\Plane{O{} O{}}{%
     \def\tmp{#1}%
     \ifx\tmp\empty%
        \mathbb{X}^2%
    \else%
        \def\tmp{#2}%
        \ifx\tmp\empty%
            \mathbb{#1}^2%
        \else%
            \mathbb{#1}^2_{\Spaceradius}%
        \fi%
    \fi%
}%
\newcommand*{\Gensin}[1]{%
    \def\symbol{\operatorname{\mathfrak{s}}}%
    \def\tmp{#1}%
    \ifx\tmp\empty%
        \symbol%
    \else%
        \symbol (#1)%
    \fi%
}%
\newcommand*{\Gencos}[1]{%
    \def\symbol{\operatorname{\mathfrak{c}}}%
    \def\tmp{#1}%
    \ifx\tmp\empty%
        \symbol%
    \else%
        \symbol (#1)%
    \fi%
}%
\newcommand*{\Coefficient}[3]{%
    t_{#1,#2;#3}%
}%
\newcommand*{\Translation}{%
    \varphi%
}%
\newcommand*{\Line}{%
    \ell%
}%
\newcommand*{\Curvature}{%
    K%
}%
\allowdisplaybreaks

\newtheorem*{theorem}{Theorem}

\title{Theorem of three squares across three geometries}
\author{Kazuhiro Ichihara}
\address{
    Department of Mathematics,
    College of Humanities and Sciences,
    Nihon University,
    3-25-40 Sakurajosui,
    Setagaya-ku,
    Tokyo 156-8550,
    Japan
}
\email{
    ichihara.kazuhiro@nihon-u.ac.jp
}
\author{
    Akira Ushijima
}
\address{
    School of Economics and Management,
    University of Hyogo,
    Hyogo 651-2197,
    Japan
}
\thanks{The second author was supported by JSPS KAKENHI Grant Number JP22K03309.}
\email{
    ushijima@em.u-hyogo.ac.jp
}
\date{
    \today
}
\subjclass[2010]{
    Primary: 51M09,
    Secondary: 51M04.
}
\keywords{
    law of cosines,
    Pythagorean theorem,
    Euclidean triangle,
    spherical triangle,
    hyperbolic triangle.
}

\begin{document}
\begin{abstract}
In Euclidean geometry, the Pythagorean theorem is presented as an equation involving three squares.
This paper explores how analogous expressions may be identified in spherical and hyperbolic geometries.
\end{abstract}
\maketitle

Let $\triangle \mathrm{ABC}$ be a triangle with internal angles $A, B, C$ at the three vertices $\mathrm{A}, \mathrm{B}, \mathrm{C}$, respectively.
For $\triangle \mathrm{ABC}$ in Euclidean, spherical, and hyperbolic geometry, the well-known formulations of the Pythagorean theorem are stated as follows, respectively, when $A$ is the right angle:
\begin{align*}
    a^2 &= b^2 + c^2 ,&
    \cos a &= \left( \cos b \right) \left( \cos c \right) ,&
    \cosh a &= \left( \cosh b \right) \left( \cosh c \right) .
\end{align*}
Here, the symbols $a$, $b$ and $c$ denote the lengths of the sides opposite the vertices $\mathrm{A}$, $\mathrm{B}$ and $\mathrm{C}$, respectively, with $a$ being the length of the hypotenuse.
Throughout this note, we shall adhere to this standard convention for naming the sides and angles of triangles.

Although the spherical and hyperbolic formulations converge to the Euclidean case for sufficiently small triangles, they differ notably in their algebraic expressions.
In Euclidean geometry, the theorem is expressed as a relationship among the squares of the three side lengths.
Therefore, in Japan, the Pythagorean theorem is commonly known as the “Three-square theorem” (\cite{JapaneseThreeSquare}) or the “Theorem of Three Squares” (\cite{NagoyaMuseum});
these are two literal translations of the same name for the theorem in Japanese, “San-Hei-Hou no Teiri”.
By contrast, the spherical and hyperbolic versions do not retain the same algebraic resemblance to the Euclidean one.
This naturally gives rise to the following question:
Is it possible to express the Pythagorean theorem in all three geometries using a unified equation involving three squares?

The differences arise from the analogous variations in the standard laws of cosines:
\begin{align*}
    \MoveEqLeft
        \textbf{The law of cosines for a Euclidean triangle:}\\
        &a^2 = b^2 + c^2 - 2 \, b \, c \cos A;\\
    \MoveEqLeft
        \textbf{The law of cosines for a spherical triangle:}\\
        &\cos a
            = \left( \cos b \right) \left( \cos c \right)
            + \left( \sin b \right) \left( \sin c \right) \cos A;\\
    \MoveEqLeft
        \textbf{The law of cosines for a hyperbolic triangle:}\\
        &\cosh a
            = \left( \cosh b \right) \left( \cosh c \right)
            - \left( \sinh b \right) \left( \sinh c \right) \cos A.
\end{align*}

The objective of this note is to introduce a Euclidean-like law of cosines for spherical and hyperbolic triangles, and to prove these within a unified framework for the three geometries.
As a corollary, our Pythagorean theorems for spherical and hyperbolic triangles are expressed in the form of an equation involving three squares.

For our unified expressions and proofs of these formulae, the following notation is adopted:
we denote by $\Plane$ one of the following geometric planes:
the Euclidean plane $\Plane[E]$, the sphere $\Plane[S]$, or the hyperbolic plane $\Plane[H]$.
For any real number $d$, representing the distance between two points on $\Plane$, we introduce generalised trigonometric functions $\Gensin{}$ and $\Gencos{}$ as follows:

Given a real constant $\Curvature$, set 
\begin{align*}
    \Gensin{d}
    &= d - \frac{\Curvature d^3}{3!} + \frac{\Curvature^2 d^5}{5!} - \cdots
    = \sum_{i=0}^{\infty} \frac{\left( - \Curvature \right)^i d^{2 i + 1}}{\left( 2 i + 1 \right) !} ,\\
    \Gencos{d}
    &= 1 - \frac{\Curvature d^2}{2!} + \frac{\Curvature^2 d^4}{4!} + \cdots
    = \sum_{i=0}^{\infty} \frac{\left( - \Curvature \right)^i d^{2 i}}{\left( 2 i \right) !}. 
\end{align*}
As can readily be observed, these are natural generalisations of the classical $\sin$ and $\cos$ functions, based on their Taylor expansions.
The constant $K$ represents the sectional curvature of the surface, although we shall not examine this in detail here.
Precisely, for $K = -1, 0, 1$, the corresponding functions are given as follows:
\begin{align*}
    \Gensin{d} \coloneqq& 
    \begin{cases*}
        d & if $\Plane = \Plane[E]$,\\
        \sin d & if $\Plane = \Plane[S]$,\\
        \sinh d & if $\Plane = \Plane[H]$;
    \end{cases*}&
    \Gencos{d} \coloneqq&
    \begin{cases*}
        1 & if $\Plane = \Plane[E]$,\\
        \cos d & if $\Plane = \Plane[S]$,\\
        \cosh d & if $\Plane = \Plane[H]$.
    \end{cases*}
\end{align*}
These functions are of independent interest; for example, they appear in the solutions to the Jacobi equation on such surfaces, as discussed in \cite[Example 4.8.4]{Elementarydifferentialgeometry}.

Using the functions above, we provide a unified description of trigonometric identities across the three geometries:

\begin{theorem}
For a triangle $\triangle \mathrm{ABC}$ on $\Plane$ 
with vertices $\mathrm{A}$, $\mathrm{B}$ and $\mathrm{C}$, the following hold. 
\begin{align*}
    \MoveEqLeft \textbf{The law of cosines:}\\
        & \Coefficient{b}{c}{a}
            \Gensin{a}^2
            = \Coefficient{c}{a}{b}
            \Gensin{b}^2
            + \Coefficient{a}{b}{c}
            \Gensin{c}^2
            - 2 \Gensin{b} \Gensin{c} \cos A;\\
        & \Coefficient{c}{a}{b}
            \Gensin{b}^2
            = \Coefficient{a}{b}{c}
            \Gensin{c}^2
            + \Coefficient{b}{c}{a}
            \Gensin{a}^2
            - 2 \Gensin{c} \Gensin{a} \cos B;\\
        & \Coefficient{a}{b}{c}
            \Gensin{c}^2
            = \Coefficient{b}{c}{a}
            \Gensin{a}^2
            + \Coefficient{c}{a}{b}
            \Gensin{b}^2
            - 2 \Gensin{a} \Gensin{b} \cos C;\\[1ex]
    \MoveEqLeft \textbf{The Pythagorean theorem (when the angle $A$ is the right angle):}\\
        & \Coefficient{b}{c}{a}
            \Gensin{a}^2
            = \Coefficient{c}{a}{b}
            \Gensin{b}^2
            + \Coefficient{a}{b}{c}
            \Gensin{c}^2 .
\end{align*}
Here the coefficients $\Coefficient{a}{b}{c}$, $\Coefficient{b}{c}{a}$
and $\Coefficient{c}{a}{b}$ are defined respectively as follows:
\begin{align*}
    \Coefficient{a}{b}{c} &\coloneqq \frac{\Gencos{a} + \Gencos{b}}{1 + \Gencos{c}} ,&
    \Coefficient{b}{c}{a} &\coloneqq \frac{\Gencos{b} + \Gencos{c}}{1 + \Gencos{a}} ,&
    \Coefficient{c}{a}{b} &\coloneqq \frac{\Gencos{c} + \Gencos{a}}{1 + \Gencos{b}} .
\end{align*}    
\end{theorem}

Our proof of the theorem follows a similar approach to that presented in \cite[Section 4.9]{Elementarydifferentialgeometry}.
See also the web page titled “Law of sines: uniform proof of Euclidean, spherical \& hyperbolic cases” on Mathematics Stack Exchange (posted 13 Dec. 2014).
There, using the same geometric setup, the standard laws of sines and cosines are derived through calculations involving the $y$- and $z$-coordinates in a unified manner.
In contrast, our proof relies primarily on calculations involving the $x$-coordinate.

\begin{proof}[Proof of Theorem]
Our Pythagorean theorem is derived directly from our law of cosines by setting the angle $A$ as the right angle. 
Thus, we will provide a concise and unified proof of it. 

Note that each geometry $\Plane$ is realized in the Euclidean 3-space $\RealSpace$ as follows:
\begin{align*}
    \Plane[E]
        &= \Set*{ (x,y,z) \in \RealSpace \given z = 1} ,\\
    \Plane[S]
        &= \Set*{ (x,y,z) \in \RealSpace \given x^2 + y^2 + z^2 = 1} ,\\
    \Plane[H]
        &= \Set*{ (x,y,z) \in \RealSpace \given x^2 + y^2 - z^2 = -1 , z > 0} .
\end{align*}
A point in $\Plane$ with polar coordinates $(r,\theta)$ has the Cartesian coordinates: 
\[(\Gensin{r} \cos \theta, \Gensin{r} \sin \theta , \Gencos{r}).\] 
In the Euclidean and spherical cases, these coordinates follow from straightforward calculations.
In the hyperbolic case, see, for example, \cite[Section 4.9]{Elementarydifferentialgeometry} or \cite[Section 6]{Reynolds}.

Let $\Line$ be the straight line on $\Plane$ defined as the intersection with the $xz$-plane.
The ``parallel translation"
in $\Plane$
along $\Line$
with a signed distance $d$
is given as follows:
\begin{equation*}
    (x,y,z) \longmapsto
    \begin{cases*}
        (x+d, \ y, \ 1) &
            if $\Plane = \Plane[E]$,\\
        (x \cos d + z \sin d , \ y , \ -x \sin d + z \cos d) &
            if $\Plane = \Plane[S]$,\\
        (x \cosh d + z \sinh d , \ y , \ x \sinh d + z \cosh d) &
            if $\Plane = \Plane[H]$.
    \end{cases*}
\end{equation*}
Note that these are expressed in a unified way as follows: 
\begin{equation} \label{eqn:translation}
    (x,y,z) \longmapsto
        ( x \Gencos{d} + z \Gensin{d} , \ y , \ -K x \Gensin{d} + z \Gencos{d} ) 
\end{equation}
More precisely, the ``parallel translation" corresponds to a rotation in the spherical case and a hyperbolic translation in the hyperbolic case, respectively. 
It can be confirmed, by direct calculation, that the parallel translation is indeed an isometry of $\Plane$ which fixes the line $\Line$.
In the cases where $\Plane=\Plane[S]$ and $\Plane=\Plane[H]$, this may be verified by using the formulae for the distance $d$ between two points $\mathrm{P} (p_1, p_2, p_3)$ and $\mathrm{Q} (q_1, q_2, q_3)$ on $\Plane[S]$ and $\Plane[H]$, which are $\cos d = p_1 q_1 + p_2 q_2 + p_3 q_3$ and $\cosh d = - p_1 q_1 - p_2 q_2 + p_3 q_3$, respectively.

Consider a triangle $\triangle \mathrm{ABC}$ in $\Plane$ with vertices $\mathrm{A}$, $\mathrm{B}$ and $\mathrm{C}$, where $\mathrm{C}$ is at $(0,0,1)$, $\mathrm{B}$ lies on $\Line$ with its polar coordinates $(a,\pi)$, and the polar coordinates of $\mathrm{A}$ are $(b , \pi - C)$. 
Let $\Translation$ be the parallel translation along $\Line$ with distance $a$.
See Figure~\ref{Fig}. 
\begin{figure}[htb]
    \begin{picture}(250,95)(0,0)
    \put(0,0){
\begin{tikzpicture}
 \coordinate (O) at (0,0); 
 \coordinate(XS)at(-4,0); \coordinate(XL)at(4,0); \draw[semithick,->,>=stealth](XS)--(XL)node[right]{$x$}; 
 \coordinate(YS)at(0,-0.5); \coordinate(YL)at(0,2.5); \draw[semithick,->,>=stealth](YS)--(YL)node[above]{$y$}; 
 
 \coordinate [label=above:A] (A) at (-1,2); 
 \coordinate [label=below:B] (B) at (-3,0); 
 \coordinate [label=below right:{$\mathrm{C} = \Translation ( \mathrm{B} )$}] (C) at (0,0); 
 \foreach \P  in {A,B,C} \fill[black] (\P) circle (0.06);  
 \draw (A)--(B)--(C)--cycle; 

 \coordinate [label=above right:{$\mathrm{A}' := \Translation ( \mathrm{A} )$}] (A') at (2,2); 
 \coordinate  (B') at (0,0); 
 \coordinate  (C') at (3,0); 
 \foreach \P  in {A',B',C'} \fill[black] (\P) circle (0.06);  
 \draw (A')--(B')--(C')--cycle; 

\draw[thick,->,>=stealth](-0.4,1.4)--(1.1,1.4);

\end{tikzpicture}
}
\end{picture}
\caption{}\label{Fig}
\end{figure}

We calculate the coordinates of
$\mathrm{A}' \coloneqq \Translation ( \mathrm{A} )$
in the following two ways:
\begin{enumerate}[label=(\roman*)]
    \item \label{itm:CoordinateA}
        Since the polar coordinates of $\mathrm{A}$ 
        are $(b , \pi - C)$,
        its Cartesian coordinates are\\[5pt]
        $
        (\Gensin{b} \cos (\pi - C) , \ 
        \Gensin{b} \sin (\pi - C) , \ 
        \Gencos{b}) 
        = 
        (- \Gensin{b} \cos C , \ 
        \Gensin{b} \sin C , \ 
        \Gencos{b})
        $.\\[5pt]
        Applying the translation formula \eqref{eqn:translation}, the $x$-coordinate of $\mathrm{A}'$ is 
        \[ \Gencos{b} \Gensin{a} - \Gensin{b} \Gencos{a} \cos C .\]
    \item \label{itm:CoordinateA'}
        Since $\Translation (\mathrm{B})$ is at $(0,0,1)$ and $\Translation (\mathrm{C})$ lies on $\Line$ with a positive $x$-coordinate, the polar coordinates of $\mathrm{A}'$ are $(c,B)$.
        Thus the $x$-coordinates of $\mathrm{A}'$ is 
        \[\Gensin{c} \cos B.\]
\end{enumerate}

From these two expressions for the $x$-coordinate of $\mathrm{A}'$, we obtain the following:
\begin{align} 
    \Gencos{b} \Gensin{a} - \Gensin{b} \Gencos{a} \cos C
        &= \Gensin{c} \cos B ,\label{eqn:xcoord}
\end{align}

An additional equation is given by calculating the $x$-coordinate of $\Translation^{-1} (\mathrm{A}') = \mathrm{A}$: 
\begin{equation*}
    - \Gencos{c} \Gensin{a} + \Gensin{c} \Gencos{a} \cos B
        = - \Gensin{b} \cos C .
\end{equation*}
The left-hand side comes from the parallel translation $\Translation^{-1}$ along $\Line$ by a distance $-a$, along with
the Cartesian coordinates of $\mathrm{A}'$ as derived in \ref{itm:CoordinateA'},
while the right-hand side comes from the expression given in \ref{itm:CoordinateA}. 
Subtracting this equation from \eqref{eqn:xcoord} gives the following:
\begin{equation}\label{1LoC}
\Coefficient{b}{c}{a} \Gensin{a}
 = \Gensin{c} \cos B + \Gensin{b} \cos C , 
\end{equation}
which is sometimes called \textit{the first law of cosines} in Euclidean geometry. 

Now, our law of cosines can be derived from Equation~\eqref{1LoC} as follows. 
In fact, by varying the vertices, it gives rise to the following system of linear equations:
\begin{equation*}
    \begin{pmatrix}
        \Gensin{b} & \Gensin{a} & 0\\
        0 & \Gensin{c} & \Gensin{b}\\
        \Gensin{c} & 0 & \Gensin{a}
    \end{pmatrix}
    \begin{pmatrix}
        \cos A\\ \cos B\\ \cos C
    \end{pmatrix}
    =
    \begin{pmatrix}
        \Coefficient{a}{b}{c} \Gensin{c}\\
        \Coefficient{b}{c}{a} \Gensin{a}\\
        \Coefficient{c}{a}{b} \Gensin{b}
    \end{pmatrix} .
\end{equation*}
Solving this system
for $\cos A$, $\cos B$ and $\cos C$
provides our law of cosines; 
\begin{align*}
\Coefficient{b}{c}{a}
            \Gensin{a}^2
            = \Coefficient{c}{a}{b}
            \Gensin{b}^2
            + \Coefficient{a}{b}{c}
            \Gensin{c}^2
            - 2 \Gensin{b} \Gensin{c} \cos A;\\
        \Coefficient{c}{a}{b}
            \Gensin{b}^2
            = \Coefficient{a}{b}{c}
            \Gensin{c}^2
            + \Coefficient{b}{c}{a}
            \Gensin{a}^2
            - 2 \Gensin{c} \Gensin{a} \cos B;\\
        \Coefficient{a}{b}{c}
            \Gensin{c}^2
            = \Coefficient{b}{c}{a}
            \Gensin{a}^2
            + \Coefficient{c}{a}{b}
            \Gensin{b}^2
            - 2 \Gensin{a} \Gensin{b} \cos C. 
\end{align*}

\vspace{-.5cm}

\end{proof}

\bibliographystyle{vancouver}

\end{document}